\documentclass[12pt]{amsart}
\usepackage{amscd,amssymb}

\theoremstyle{definition}

\theoremstyle{remark}

\numberwithin{equation}{section}

\newcommand{\Z}{{\mathbf Z}}

\newcommand{\R}{{\mathbf R}}
\newcommand{\C}{{\mathbf C}}

\newcommand{\RR}{{\mathcal {R}}}

\renewcommand{\L}{{\mathcal {L}}}

\DeclareMathOperator{\rank}{{rank}}

\DeclareMathOperator{\sign}{{sign}}
\DeclareMathOperator{\End}{{End}}

\DeclareMathOperator{\me}{{\operatorname{Mon_E}}}

\DeclareMathOperator{\px}{{\widehat{\Z\pi_\xi}}}
\DeclareMathOperator{\sx}{{\Sigma_\xi^{-1}(\Z\pi)}}

\DeclareMathOperator{\nv}{{{\operatorname{\mathbf {Nov}}}}}

\DeclareMathOperator{\kk}{{\mathbf {k}}}
\DeclareMathOperator{\V}{{\mathcal V}}

\DeclareMathOperator{\Hom}{{\operatorname{Hom}}}

\DeclareMathOperator{\GL}{{\operatorname{GL}}}
\DeclareMathOperator{\ml}{{\operatorname{Mon_L}}}
\DeclareMathOperator{\In}{{\operatorname{{Int}}}}

\begin{document}

\title[Morse - Novikov critical point theory]
{Morse - Novikov critical point theory, \\
Cohn localization and Dirichlet units}

\author[M.~Farber]{M.~Farber}
\address{Department of Mathematics, Tel Aviv University, Tel Aviv, 69978, Israel}
\email{Farber@math.tau.AC.IL}
\date{\today}

\subjclass{Primary 57Q10;  Secondary 53C99}
\keywords{Morse theory, closed 1-forms, Novikov inequalities, 
Cohn localization}
\thanks{The research was partially supported by a grant from the 
US - Israel Binational Science Foundation and by
the Herman Minkowski Center for Geometry; 
a part of this work was done while the author
visited Max-Planck Institut f\"ur Mathematik in Bonn}

\begin{abstract}
In this paper we construct a Universal chain complex, 
counting zeros of closed 1-forms on a manifold. The Universal complex is a refinement of
the well known Novikov complex;
it relates the homotopy type of the manifold, after a suitable 
noncommutative localization, 
with the numbers of zeros of different indices 
which may have closed 1-forms within a given cohomology class. 
The Main Theorem of the paper generalizes the result of a joint paper with A. Ranicki 
\cite{FR}, which treats the
special case of closed 1-forms having integral cohomology classes. 
The present paper also describes a number of new inequalities,
giving topological lower bounds on the minimum number of zeros of closed 
1-forms. 
In particular, such estimates are provided by
the homology of 
flat line bundles with monodromy described by complex numbers, 
which are not Dirichlet units. 
\end{abstract}

\maketitle

\section{\bf The Main Theorem}

\subsection{Basic definitions}
Let $M$ be a smooth manifold, and let $\omega$ be a closed 1-form on $M$, $d\omega=0$.
A point $p\in M$ is {\it a zero} of $\omega$ if $\omega$ vanishes at this point, i.e. $\omega_p=0$. 

A zero $p\in M$ is called {\it nondegenerate} if $\omega$, 
viewed as a map $M\to T^\ast M$, is transversal to the zero section $M\subset T^\ast M$ of the 
cotangent bundle. As is well-known, this condition is equivalent to the requirement that 
in a neighborhood $U$ of $p$ we may write $\omega=df$, where $f:U\to \R$
is a smooth function and $p$ is a non-degenerated (Morse) critical point of $f$. 
The {\it Morse index }
of $p$ is well defined (as the Morse index of $p$ as the critical point of $f$). 

A closed 1-form $\omega$ is called {\it Morse} if
all its zeros are non-degenerate.

\subsection{Morse theory} 
Let $\pi$ denote the fundamental group $\pi_1(M)$. Fix a cell decomposition of $M$ and let $C=C_\ast(\tilde M)$ denote the free
$\Z\pi$-chain complex of cellular chains in the universal covering $\tilde M$. 

From the classical Morse - Smale - Thom theory we know that any
Morse function $f: M\to \R$ determines (via the trajectories of the gradient flow) 
a decomposition of $M$ into open
disks (one for each critical point of $f$); this yields a free finitely generated chain complex
$C^f$ of $\Z\pi$-modules with the following properties:
{\it 
(a) $C^f$ is chain homotopy equivalent to $C_\ast(\tilde M)$;
(b) each chain module $C^f_j$ has a free basis which is in a canonical 
one-to-one correspondence with the 
critical points of $f$ having Morse index $j$, where $j=0, 1, \dots, n=\dim M$.}

The main purpose of this paper is to construct an analog $C^\omega$ of the chain complex $C^f$ 
for closed 1-forms $\omega$; if the form $\omega$ is exact $\omega =df$ then our complex $C^\omega$ coincides with $C^f$.

Namely, for any closed 1-form $\omega$ with Morse zeros on $M$ we will construct 
a chain complex $C^\omega$ with the following properties:

{\it (i) $C^\omega$ is a chain complex of free modules over a localization of the group ring $\Z\pi$;  

(ii) $C^\omega$ is chain homotopy equivalent to the localized chain complex $C_\ast(\tilde M)$; 

(iii) there is a canonical free basis of $C^\omega_j$, which is in a one-to-one 
correspondence with the zeros of the form $\omega$, having index $j$, where $0\le j\le n=\dim M.$}

Note that the localization, which is mentioned in (i) is a {\it non-commutative localization}; 
the theory of such localization was developed by P.M. Cohn \cite{C}. Roughly, 
it consists in inverting
a class of square matrices, and not single elements, as in the commutative case. 

The class of square matrices, which we invert, 
depends on the cohomology class of $\xi=[\omega]\in H^1(M;\R)$ of the closed 1-form $\omega$,
cf. below.

In order to state our main result, we will make the following definitions.

\subsection{$\xi$-negative matrices} We will view the cohomology class $\xi\in H^1(M;\R)$ as a
homomorphism $\xi: \pi\to \R$. We have $\xi(gg') =\xi(g) +\xi(g')$ for $g,g'\in \pi$.
An element $\alpha\in \Z\pi$ will be called {\it $\xi$-negative} if 
$\alpha=\sum n_j g_j$ (finite sum) where $n_j\in \Z$ and $\xi(g_j)<0$ for all $j$. 
An $m\times m$-matrix $A$ over the
group ring $\Z\pi$ will be call {\it $\xi$-negative} if all its entries are $\xi$-negative.

\subsection{\bf Main Theorem (the first form)} {\it Let $M$ be a closed smooth manifold with 
$\pi=\pi_1(M)$,
and let $\xi\in H^1(M;\R)$ be a cohomology class.
Let $\rho: \Z\pi\to \RR$ be a ring homomorphism with the following
property: for any $\xi$-negative square matrix $A$ over $\Z\pi$ the matrix $\rho(I+A)$ over $\RR$
is invertible. (Here $I$ denotes the unit matrix of the appropriate size). 
Then for any closed 1-form
$\omega$ on $M$ having only Morse zeros and representing the class $\xi$,
there exists a free chain complex $C^{\omega,\rho}$
of $\RR$-modules, such that $C^{\omega,\rho}$ is chain homotopy equivalent to 
\begin{eqnarray}
\RR\otimes_{\Z\pi}C_\ast(\tilde M)
\end{eqnarray}
and each $\RR$-module
$C^{\omega,\rho}_j$ has a canonical free basis, which is in a one-to-one correspondence 
with the zeros of the form $\omega$
having index $j$, where $0\le j\le n=\dim M$.}

Here $C_\ast(\tilde M)$ denotes the
chain complex of the 
universal covering of $M$ corresponding to a cell decomposition of $M$, and $\RR$ is viewed as
a right $\Z\pi$-module via homomorphism $\rho$.

Note that in the case $\xi=0$ (the usual Morse theory of functions) 
there are no $\xi$-negative matrices and hence
any ring homomorphism $\rho: \Z\pi\to \RR$ satisfies the condition of Theorem 1.4.

In the special case of closed 1-forms having integral cohomology classes 
(the circle-valued Morse theory)
the Main Theorem was established in a joint paper \cite{FR} with A.A. Ranicki.

We will show in sections \S 3, 4, 5 that the Main Theorem implies new 
Novikov type inequalities. The Novikov inequalities also follow immediately, cf. \S 2.

\subsection{ Example: The Novikov-Sikorav completion} An example of a ring 
homomorphism 
$\rho: \Z\pi\to \RR$ with the property that for any $\xi$-negative square matrix $A$ over $\Z\pi$
the matrix $\rho(I+A)$ with entries in $\RR$ is invertible, is provided by the homomorphism
\begin{eqnarray}
\rho: \Z\pi\to \px,
\end{eqnarray}
known as the {\it Novikov-Sikorav completion},
which we now recall. It was first introduced by S.P. Novikov \cite{N2}
for free abelian quotients
and by J.-C. Sikorav \cite{S1} in the general case.

Elements of the ring $\px$ are represented by formal sums, possibly infinite,
$\alpha = \sum n_i g_i$, where $n_i\in \Z$ and $g_i\in \pi$, satisfying the following condition:
{\it for any $c\in \R$ the set $\{i; \, \xi(g_i)\ge c\}$ is finite. } The addition and the 
multiplication are given
by the usual formulae; for example, the product of $\alpha=\sum n_i g_i\in \px$ and 
$\beta=\sum m_j h_j\in \px$ is given by
\begin{eqnarray}
\alpha\cdot\beta = \sum_{i,j} (n_i m_j) (g_i h_j).
\end{eqnarray}
The ring homomorphism $\rho: \Z\pi\to \px$ is the inclusion. If $A$ is a $\xi$-negative square matrix 
over the ring $\Z\pi$, then the power series 
$(I+A)^{-1} = I - A + A^2 - \dots$
converges in $\px$
and hence the matrix $I+A$ is invertible in $\px$. Thus the Novikov-Sikorav completion
satisfies the condition of Theorem 1.4.

\subsection{Example: the Novikov ring} Consider the ring $\nv$ consisting of formal power series
\[\sum_{j\ge 1} n_j t^{\gamma_j},\]
where the coefficients $n_j\in \Z$ are integers and the exponents $\gamma_j\in \R$ are real and satisfy
the condition $\lim_{n\to \infty}\gamma_n = -\infty$. 

Any cohomology class $\xi:\pi \to \R$
gives a representation $\rho_\xi : \Z\pi \to \nv$, which sends a group element $g\in \pi$ to the monomial
$t^{\xi(g)}\in \nv$. If $A$ is a $\xi$-negative square matrix over $\Z\pi$ then the determinant 
of the matrix $\rho_\xi (I+A)$ has the top coefficient 1 and hence $\rho_\xi (I+A)$ is invertible in 
the Novikov ring $\nv$. This shows that the homomorphism $\rho_\xi : \Z\pi \to \nv$ 
satisfies Theorem 1.4.

\subsection{Example: a ring of rational functions \cite{F1}} 
Let $\RR$ be the following ring 
of rational functions in the indeterminate $t$.
Elements of $\RR$ are rational functions of the form $p(t)/q(t)$, where 
$p(t), q(t)\in Z[t,t^{-1}]$ are Laurent
polynomial with integral coefficients and the denominator $q(t)$ has top coefficient 1, i.e.
$q(t)= t^m+b_1t^{m-1}+\dots +b_m$, where $b_i\in \Z$. This ring was introduced in \cite{F1}.
It is shown in \cite{F1} that $\RR$ is a principle ideal domain.

Suppose that the cohomology class $\xi: \pi\to \Z\subset \R$ is integral.
It defines the representation $\rho_\xi: \Z\pi\to \RR$, sending a group element
$g\in \pi$ to $t^{\xi(g)}\in\RR$.
The homomorphism $\rho_\xi$ 
satisfies the condition of Theorem 1.4, since for any $\xi$-negative square matrix $A$ over
$\Z\pi$ the matrix $\rho_\xi(I+A)$ will be an $\RR$-matrix, having a Laurent polynomial
with integer coefficients and the top coefficient 1 as its determinant.

\subsection{Example: algebraic integers} Suppose,
that the cohomology class $\xi$ is integral. 
Given a complex number $a\in \C^\ast$, we obtain a ring homomorphism
$\rho_a: \Z\pi \to \C$ given by 
\begin{eqnarray}
\rho_a(\alpha) = \sum n_j a^{\xi(g_j)}\in \C, \qquad\text{where}\quad \alpha = \sum n_j g_j\in \Z\pi.
\end{eqnarray}
Suppose that {\it the number $a$ is not
an algebraic integer.} Then {\it the homomorphism $\rho_a$ satisfies the condition of Theorem 1.4} 
i.e. for any $\xi$-negative square matrix $A$ over $\Z\pi$
the matrix $\rho_a(I+A)$ with entries in $\C$ is invertible. Indeed, the determinant of the matrix 
$\rho_a(I+A)$ clearly can be represented as 
an evaluation of a polynomial 
\[p(\lambda)= 1+\beta_1\lambda^{-1}+\beta_2\lambda^{-2}\dots \beta_k\lambda^{-k}, 
\quad \beta_j\in \Z\]
at the value of the parameter
$\lambda = a$ and hence it is nonzero, and the matrix is invertible.

\subsection{Example: the field of rational functions} Here is an example of a ring homomorphism
which satisfies the condition of Theorem 1.4 for any class $\xi\in H^1(M;\R)$.

Let $\xi_1, \dots, \xi_r: \pi \to \Z$ form a basis of the free abelian group $\Hom(\pi;\Z)$.
Let $\kk$ be a field and let $\kk(t_1, \dots, t_r)$ denote the field of rational functions
in the indeterminates $t_1, \dots, t_r$. We obtain a ring homomorphism 
$\rho: \Z\pi\to \kk(t_1, \dots, t_r)$ which sends a group element $g\in \pi$ to the monomial
$t_1^{\xi_1(g)}\cdot t_2^{\xi_2(g)}\cdot \dots \cdot t_r^{\xi_r(g)}$. 

Any cohomology class
$\xi:\pi \to \R$ is represented by a sum 
$\xi =\sum_{i=1}^r \lambda_i\xi_i,$
where the coefficients $\lambda_i\in \R$ are real. {\it A $\xi$-weight} of a monomial $t_1^{m_1}\dots t_r^{m_r}$
is defined as the number 
$m_1\lambda_1 + \dots +m_r\lambda_r\, \in \, \R.$
If $A$ is a $\xi$-negative 
square matrix over the group ring $\Z\pi$ then the determinant 
$\det(\rho(I+A))\in \kk(t_1, \dots, t_r)$ is nonzero,
since it is represented by a Laurent
polynomial of the form $1+$ terms having negative $\xi$-weight. Hence, 
the matrix $\rho(I+A)$ is invertible.

\subsection{The Cohn localization}
We may restate our main theorem 1.4 in a different form using the Cohn localization \cite{C}. 
Given a group $\pi$ and a homomorphism $\xi:\pi \to \R$, where $\xi(gg') = \xi(g)+\xi(g')$ for
$g, g'\in \pi$, consider {\it the universal Cohn localization} 
\begin{eqnarray}\rho_\xi : \Z\pi \to \sx,\end{eqnarray}
inverting the class $\Sigma_\xi$ 
of square matrices of the form $I+A$, where $A$ is a $\xi$-negative square matrix.
Recall that this means that the homomorphism (1.5) satisfies the following property: firstly, 
any matrix $\rho_\xi(I+A)$, where $A$ is $\xi$-negative, is invertible over $\sx$, and, 
secondly, it is a {\it universal} homomorphism having this property, i.e. for any ring
homomorphism $\rho: \Z\pi\to \RR$, 
inverting matrices of the form $I+A$, where $A$ is $\xi$-negative,
there exists a unique ring homomorphism $\phi: \sx \to \RR$ such that the following diagram commutes.

\begin{figure}[h]
  \begin{center}
\begin{picture}(9,2.8)
\linethickness{0.1mm}
\put(2.7,0.5){$\Z\pi$}
\put(4,2.){$\sx$}
\put(6,0.5){$\RR$}
\put(3.7,0.6){\vector(1,0){1.9}}
\put(4.4,0.1){$\rho$}
\put(3.3,1.0){\vector(1,1){0.7}}
\put(3.3,1.5){$\rho_\xi$}
\put(5,1.7){\vector(1,-1){0.7}}
\put(5.5, 1.4){$\phi$}
\end{picture}
\end{center}
\end{figure}

We refer to the book of P.M. Cohn \cite{C}, where existence of localization (1.5) is proven.
In particular, there is a canonical ring homomorphism 
\begin{eqnarray}
\sx\to \px,\end{eqnarray}
extending the inclusion $\Z\pi \to \px$. This implies that {\it the homomorphism (1.5) is injective.}
However, it is not known when (1.6) is injective.
Intuitively, the image of (1.6) consists of {\it "rational power series"}. This is a consequence 
of a characterization of elements of the Cohn localization as components of solutions of linear
system of equations, cf. \cite{C}, chapter 7.

\subsection{Main Theorem the second form)} {\it
Let $M$ be a closed smooth manifold with 
$\pi=\pi_1(M)$,
and let $\xi\in H^1(M;\R)$ be a cohomology class.
Let $\rho_\xi: \Z\pi\to \sx$ be the Cohn localization of the group ring determined by the class $\xi$
(cf. above). Then for any closed 1-form
$\omega$ on $M$ having only Morse zeros and representing the class $\xi$,
there exists a free chain complex $C^\omega$ consisting 
of $\sx$-modules, such that $C^\omega$ is chain homotopy equivalent to the localized complex \begin{eqnarray}
\sx\otimes_{\Z\pi}C_\ast(\tilde M)
\end{eqnarray}
and each $\sx$-module
$C^\omega_j$, where $0\le j\le n=\dim M$, 
has a canonical free basis, which is in a one-to-one correspondence 
with the zeros of the closed 1-form $\omega$
having index $j$.}

Note that in the case $\xi=0$ (when one studies Morse functions) the set of $\xi$-negative matrices
is empty and so the Cohn localization (1.5) is just the identity map.

It is clear that Theorem 1.11 is equivalent to Theorem 1.4.

A proof of Theorem 1.11 will be given in \S 8.

\subsection{The Universal Complex} 
We will call the chain complex $C^\omega$ of $\sx$-modules,
which appears in Theorem 1.11, 
{\it the universal complex determined by the closed 1-form $\omega$.} 

We will restate Main Theorem 1.11 in the form useful for some applications.

\subsection{Main Theorem (the third form)} {\it
Let $M$ be a closed smooth manifold and let $\xi\in H^1(M;\R)$ be a cohomology class. 
Let $\RR$ be a ring and let $X$ be a 
$(\RR, \Z\pi)$-bimodule (where $\pi=\pi_1(M)$), satisfying the following condition:
for any $\xi$-negative $m\times m$ matrix $A$ over 
$\Z\pi$ the map 
\begin{eqnarray}
X^m \to X^m,\quad \text{where}\quad 
(x_1, \hdots, x_m)\mapsto 
(x_1, \hdots, x_m)\cdot (I+A),
\end{eqnarray}
is an $\RR$-isomorphism. Then for any closed 1-form
$\omega$ on $M$ having only Morse zeros and representing the class $\xi$,
there exists a chain complex $C^{\omega, X}$ of $\RR$-modules, which is homotopy equivalent
to $X\otimes_{\Z\pi}C_\ast(\tilde M)$ and
such that each $C^{\omega, X}_j$ is isomorphic to a direct sum of $c_j(\omega)$ copies of $X$. 
Here $c_j(\omega)$ denotes the number of zeros of $\omega$
having Morse index $j$, where $j=0, 1, 2, \dots, n$.}

Theorem 1.13 obviously follows from Theorem 1.11. The condition on $X$ implies that its
$(\RR,\Z\pi)$-bimodule structure can be extended to a $(\RR, \sx)$-bimodule structure. Hence
we may set $C^{\omega, X} = X\otimes_{\sx} C^\omega$, where $C^\omega$ is the Universal 
Complex 1.12.

\subsection{The Novikov complex} 
S. P. Novikov conjectured in \cite{N1}, \cite{N2}, that any closed 1-form 
must generate a chain complex over the Novikov ring $\nv$, which has the zeros 
of the closed 1-form as the free basis of the module of the chains and has prescribed
chain homotopy type. 

Before the present paper, 
existence of the Novikov complex was rigorously proven only for closed
1-forms with integral cohomology classes, see \cite{P1}.

{\it Theorem 1.4 and example 1.5 give a proof of existence 
of the Novikov complex in full generality, without any restrictions on the cohomology class. }

The advantage of the Universal complex $C^\omega$, compared to the Novikov complex, is that
$C^\omega$ produces many different {\it "Novikov complexes"} 
$C^{\omega, \rho}$, one for each $\Sigma_\xi$-inverting representation $\rho$, cf. Theorem 1.4;
the original Novikov complex corresponds to a specific $\rho$, described in 1.5.

\subsection{} A. Pajitnov announced  in \cite{P2} a theorem that the incidence coefficients of the
Novikov complex of a Morse form with integral cohomology class
belong to the image of the Cohn localization in the Novikov-Sikorav completion; preprint
\cite{P3} contains a proof. 

The Theorem of \cite{P2}, \cite{P3}
provides an interesting information about the {\it dynamical properties 
of the gradient vector fields of closed 1-forms}. 
However, this result does not imply existence of a lift of the Novikov complex
to a complex over the Cohn localization $\sx$, even in the studied in \cite{P3}
case of forms with integral cohomology classes.

\subsection{} The plan of the paper is as follows. In \S 2 we show how the classical Novikov
inequalities follow from our Main Theorem. In sections \S\S 3, 4, 5 we prove some 
new theorems giving
topological lower bounds on the Morse numbers of closed 1-forms. 
In \S 6 we briefly remind the
basic facts of the Morse theory for manifolds with corners, which will be later 
used in the proof of the Main Theorem. In \S 7 we describe the chain collapse technique, developed
in \cite{FR}, which we also need as an important component of the proof. The last section \S 8
contains the proof of the Main Theorem 1.11.

\subsection{}  I would like to thank the Max-Planck Institut f\"ur Mathematik in Bonn for hospitality.

\section{\bf The Novikov inequalities}

In this section we illustrate our Main Theorem by showing how the well-known Novikov
inequalities \cite{N1}, \cite{N2} can be obtained from it. 

In the next sections we will apply the Main 
Theorem to obtain some new inequalities.

\subsection{The Novikov numbers}
Consider the Novikov ring $\nv$ and the canonical homomorphism 
$\rho_\xi : \Z\pi \to \nv$ determined
by a cohomology class $\xi: \pi \to \R$. Note that the ring $\nv$ is a principal ideal domain, cf. \cite{HS}. Given a manifold $M$ with $\pi_1(M)=\pi$ and a cell decomposition of $M$,
consider the chain complex 
\begin{eqnarray}
\nv\otimes_{\Z\pi}C_\ast(\tilde M).\end{eqnarray}
It is a complex of finitely
generated modules over $\nv$. The homology 
\begin{eqnarray}H_i(\nv\otimes_{\Z\pi}C_\ast(\tilde M))\end{eqnarray}
 is a
finitely generated module over $\nv$, 
and so it can be represented as a direct sum of a free module and a 
torsion submodule. {\it The Novikov number} 
$b_i(\xi)$ is defined as the minimal number of generators
of the free part of (2.2); {\it the Novikov number} 
$q_i(\xi)$ is defined as the minimal number of generators
of the torsion part of (2.2). 

\subsection{Corollary (the Novikov inequalities)} {\it Let $\omega$ be a closed 1-form 
with Morse zeros on a smooth
closed manifold $M$. Then the number $c_j(\omega)$ of zeros of $\omega$ having Morse index $j$
satisfies
\begin{eqnarray}
c_j(\omega) \ge b_j(\xi) + q_j(\xi) + q_{j-1}(\xi),\end{eqnarray}
where $\xi\in H^1(M;\R)$ is the cohomology class of the form $\omega$.}

\begin{proof} 
By Theorem 1.4, the complex (2.1) is chain homotopy equivalent to a free chain complex
$C^{\omega, \rho_\xi}$ over the ring $\nv$, such that
the rank of each module $C^{\omega, \rho_\xi}_j$
equals $c_j(\omega)$ for any $j$. Hence the Novikov numbers $b_i(\xi)$ and $q_i(\xi)$ can be 
computed starting from the complex $C^{\omega, \rho_\xi}$. 
The inequality (2.3) is now standard for complexes over principal ideal domains.
This completes the proof. 
\end{proof}

\section{\bf Dirichlet units and inequalities for critical points}

In this section we describe some new inequalities for critical points of closed 1-forms. They are 
easy consequences of our Main Theorem, cf. 1.4, 1.11 and 1.13. 

\subsection{Dirichlet units} Recall that a complex number $a\in \C^\ast$ is called {\it a Dirichlet
unit} if it is a unit of the ring of integers of an algebraic number field. Equivalently, a Dirichlet unit
is an algebraic integer such that its inverse $a^{-1}$ is also an algebraic integer. 
Dirichlet unit $a\in \C$ is a root of a monic polynomial
\[
a^k + \beta_1a^{k-1} + \beta_2a^{k-2} +\dots + \beta_k\, =\, 0
\] 
with $\beta_1, \beta_2, \dots, \beta_{k-1}\in \Z$ and $\beta_k=\pm 1$; this property clearly characterizes the Dirichlet units.

In \cite{F2} it is shown that Dirichlet units play a crucial role in the critical point theory
of closed 1-forms, when we do not assume any non-degeneracy 
({\it Lusternik - Schnirelman type theory}).

\subsection{Notation} 
We will consider flat complex vector bundles $E$ over a closed manifold $M$. 
As is well known, a flat vector bundle is
determined by its monodromy, a linear representation of the fundamental group $\pi_1(M,x_0)$ 
on the fiber $E_0$ over the base point $x_0$, which is given by the parallel transport along loops.
For example, a flat line bundle is determined by a homomorphism $H_1(M;\Z)\to \C^\ast$, where $\C^\ast$ is considered 
as a multiplicative abelian group.

{\it A lattice} $\L\subset V$ in a finite dimensional vector space $V$ is a finitely 
generated subgroup with $\rank \L =\dim_\C V$.
We will say that a complex flat bundle $E\to M$ of rank $m$ {\it admits an integral lattice} if its monodromy representation $\pi_1(M,x_0)\to \GL_{\C}(E_0)$ is conjugate to a homomorphism 
$\pi_1(M,x_0)\to \GL_{\Z}(\L_0)$, where $\L_0\subset E_0$ is a lattice in the fiber. This condition
is equivalent to the assumption that $E$ is obtained from a local system of finitely generated free 
abelian groups over $M$ by tensoring on $\C$.

Let $\xi\in H^1(M;\Z)$ be an integral cohomology class. Given a complex number $a\in \C^\ast$,
we will consider the complex flat line bundle over $M$ with the following property: the monodromy
along any loop $\gamma\in \pi_1(M)$ is the multiplication by $a^{\langle\xi,\gamma\rangle}$. We will denote 
this bundle by $a^\xi$.

\subsection{Theorem} {\it Let $M$ be a closed smooth manifold and let $\xi\in H^1(M;\Z)$ 
be an integral cohomology class. Let $E\to M$ be a complex flat bundle admitting an integral lattice.
Let $a\in \C^\ast$ be a complex number, which is not a Dirichlet unit.
Then for any closed 1-form $\omega$ on $M$ having 
Morse type zeros and lying in the class $\xi$, the number $c_p(\omega)$ 
of zeros of $\omega$ having index $p$ satisfies
\begin{eqnarray}
c_p(\omega)\, \ge\,  \frac{\dim_\C H_p(M;a^\xi\otimes E)}{\dim E},\qquad p=0,1,2, \dots.\end{eqnarray} 
Moreover, 
\begin{eqnarray}\qquad\sum_{j=0}^p (-1)^j c_{p-j}(\omega) \, \ge\,  
\sum_{j=0}^p (-1)^j \frac{\dim_\C H_{p-j}(M;a^\xi\otimes E)}{\dim E},\,\,\, p=0,1,2, \dots.
\end{eqnarray}}

Theorem 3.3 gives interesting estimates already in the simplest case when $E$ is taken to be 
the trivial flat line bundle.

\begin{proof} Suppose first that $a\in \C$ is not an algebraic integer. 
Let $\rho_E: \pi_1(M)\to \GL_{\C}(E_0)$ be the monodromy representation 
of the flat bundle $E\to M$,
where $E_0$ is the fiber over the base point. Note that holds
$\rho_E(gg') = \rho_E(g')\circ\rho_E(g),\quad g,g'\in \pi,$
i.e. $\rho_E$ defines a {\it right} action of $\pi$ on $E_0$.
Let $\rho_a:\Z\pi \to \C$ be the representation defined
in section 1.8. We obtain a ring homomorphism 
$\rho= \rho_a\otimes \rho_E:\Z\pi \to \End(E_0),$
and $E_0$ becomes a $(\C, \Z\pi)$-bimodule via $\rho$, satisfying the conditions of Theorem 1.13 
(arguments proving that are the same as in 1.8). Hence by Theorem 1.13 for any closed
1-form $\omega$ on $M$ lying in class $\xi$ and having only Morse zeros, 
there exists a chain complex 
$C_\ast$ of complex vector spaces having $\dim(E)\cdot c_j(\omega)$ generators in any dimension $j$ and  computing the homology 
$H_\ast(M;a^\xi\otimes E)$. The inequalities (3.1) and (3.2) follow from existence of $C_\ast$
via the standard well-known argument, cf. \cite{M}.

Suppose now that $a^{-1}\in \C$ is not an algebraic integer. Consider the 
dual vector bundle $E^\ast\to M$. The above arguments
applied to the form $-\omega$ prove the inequalities
\[
c_p(-\omega)  \ge 
\frac{\dim H_p(M;a^{-\xi}\otimes E^\ast\otimes {\mathfrak o}_M)}{\dim E},\]
where $n=\dim M$ and 
$\mathfrak o_M$ is the orientation bundle of $M$ 
(i.e. a flat line bundle such the monodromy along any loop equals $\pm 1$ 
depending on whether the orientation of $M$ is preserved or reversed along the loop). 
Note, that both flat bundles $E^\ast$ and 
$\mathfrak o_M$ have integral lattices. Inequality (3.1) follows now using 
$c_p(-\omega) = c_{n-p}(\omega)$ and the Poincar\'e duality
$H_p(M;a^{-\xi}\otimes E^\ast\otimes {\mathfrak o}_M)\simeq H_{n-p}(M;a^\xi\otimes E)$.

Inequality (3.2) is obtained similarly. By the above arguments applied to the form $-\omega$ and
the flat vector bundle $E^\ast\otimes \mathfrak o_M$ we obtain the inequalities
\[\qquad\quad
\sum_{j=0}^{n-p} (-1)^j c_{p+j}(-\omega) \, \ge\,  
\sum_{j=0}^{n-p} (-1)^j 
\frac{\dim_\C H_{p+j}(M;a^{-\xi}\otimes E^\ast\otimes \mathfrak o_M)}{\dim E},\]
for $ p=0,1,2, \dots$
(a variant of the usual inequality with increasing indices)
and then one uses the Poincar\'e duality and $c_p(-\omega) = c_{n-p}(\omega)$ to obtain (3.2).

\end{proof}

\subsection{Remark} It is easy to show that the Betti
number $\dim_\C H_i(M; a^\xi)$ for transcendental $a\in \C$ 
equals the Novikov number $b_i(\xi)$ (and in, particular, 
it is the same for all transcendental $a$. 

\subsection{Remark} Consider the function
$$a\in\C^\ast \mapsto \dim_\C H_i(M; a^\xi\otimes E).$$
Then there exist only finitely many numbers 
$a_1, a_2, \dots, a_k \in \C^\ast$ (they are called {\it jump points}) so that the corresponding 
Betti number $\dim_\C H_i(M; a^\xi\otimes E)$ is the same for any $a\in C^\ast$ which is not 
one of the jump points.
Following \cite{BF}, let us denote by $b_i(\xi; E)$ the value of 
$\dim_\C H_i(M; a^\xi\otimes E)$ for $a$ not a jump point.
The number $b_i(\xi; E)$ is a generalization of the Novikov number $b_i(\xi)$. 
For any of the jump points $a_j$ actually holds
$$\dim_\C H_i(M; a_j^\xi\otimes E) > b_i(\xi; E),$$
i.e. {\it the jumps are always positive.} 

Suppose that the flat bundle $E$ admits an integral lattice.
Then the jump points $a_1, a_2, \dots, a_k $ are 
{\it algebraic numbers} (not necessarily algebraic integers).
If a jump point happens at a point, which is not a Dirichlet unit, then Theorem 3.1 applies and 
we obtain estimate (3.1) which is stronger 
than the inequality 
$c_i(\omega) \ge b_i(\xi; E)/\dim E.$

\subsection{Remark} {\it The inequality (3.1) is false if $a$ is a Dirichlet unit. }
To explain this, note that any Dirichlet unit $a\in \C$ is an eigenvalue of an integral square 
matrix $B=(b_{ij})$ with $\det(B)=1$. We may find a diffeomorphism of a compact smooth 
manifold $h: F\to F$ so that $h$ induces the matrix $B$ on homology 
of some dimension $k$.
Consider the mapping torus $M$, which obtained from $F\times [0,1]$ by identifying
any point $(x,0)$ with $(h(x),1)$. The manifold $M$ is naturally a smooth fiber bundle over the
circle and so it admits a closed 1-form $\omega$ with no critical points, $c_i(\omega)=0$ for all $i$.
The homology $H_\ast(M; a^\xi)$ is nontrivial if and only if the number $a$ is an eigenvalue
of the monodromy $h_\ast: H_\ast(F;\C)\to H_\ast(F;\C)$. Here $\xi$ denotes the cohomology class
of $\omega$. Hence, if $a$ is not a Dirichlet unit, we may construct $M$ so that
$H_\ast(M;a^\xi)\ne 0$ and class
$\xi$ may be realized by a closed 1-form with no critical points.

\subsection{Remark} In \cite{F4} a different proof of 
Theorem 3.3 was suggested. \cite{F4} contains examples showing that
Theorem 3.3 may produce stonger estimates than the Novikov inequalities.
\cite{F4} also contains a generalization of Theorem 3.3 for closed 1-forms
with non-isolated zeros.

\section{\bf Line bundles and Dirichlet units}

In this section we state generalization of the results of \S 3 for cohomology classes $\xi$
of higher rank. 

\subsection{} Let $M$ be a manifold. We will denote by $H$ the first homology group $H_1(M;\Z)$.
Let $\xi\in H^1(M,\R)$ be a real cohomology class. It can be viewed as a homomorphism
$\xi: H_1(M;\Z)=H\to \R$; we denote by $\ker(\xi)$ the kernel. 
Given a polynomial $p\in \Z[H]$, one defines two numbers
$d_\xi(p)$ ({\it the $\xi$-degree of $p$}) and $v_\xi(p)$ 
({\it the $\xi$-top coefficient}) as follows. Let $p=\sum_{j=1}^n \beta_j h_j$, where 
$\beta_j\in \Z$ and $h_j\in H$. Then  
$d_\xi(p)$ is defined as the maximal number $d=d_\xi(p)\in \R$ such that the sum
$v_\xi(p) = \sum \beta_j$, taken over all $j$ with $\langle\xi, h_j\rangle=d$, is nonzero.

Let $L\to M$ be a complex flat line bundle. We will assume that the monodromy of $L$ is trivial
along any loop in $M$ representing a homology class in $\ker(\xi)$. 
$L$ determines the {\it monodromy homomorphism}
\[\ml : \Z[H]\to \C,\]
which is a ring homomorphism extending the map assigning to any $h\in H$ the monodromy of $L$
along $h$. We will denote by $\mathcal I_L\subset  \Z[H]$ the kernel of the homomorphism $\ml$.

\subsection{Definition} (A) {\it
We will say that a flat complex line bundle $L\to M$ is a {\it $\xi$-algebraic integer}
if (i) the monodromy of $L$ is trivial
along any loop in $M$ representing a homology class in $\ker(\xi)$; and (ii) the ideal $\mathcal I_L$
contains a polynomial $p\in \mathcal I_L$ with $v_\xi(p)=\pm 1$.}

(B) {\it We will say that a complex flat line bundle $L\to M$ is a $\xi$-Dirichlet unit if $L$ and the
dual flat line bundle $L^\ast$ are
$\xi$-algebraic integers. } 

Now we may state a generalization of Theorem 3.3 for classes of rank $>1$:

\subsection{Theorem} {\it Let $M$ be a closed smooth manifold and let $\xi\in H^1(M;\R)$ 
be a real cohomology class. Let $E\to M$ be a flat complex vector bundle admitting an integral lattice.
Let $L\to M$ be a flat complex line bundle, which is not a $\xi$-Dirichlet unit.
Then for any closed 1-form $\omega$ on $M$ having 
Morse zeros and lying in the class $\xi$, the number $c_p(\omega)$ of zeros
of $\omega$ having index $p$ satisfies
\begin{eqnarray}
c_p(\omega)\, \ge\,  \frac{\dim_\C H_p(M;L\otimes E)}{\dim E},\qquad p=0,1,2, \dots.
\end{eqnarray} 
Moreover, 
\begin{eqnarray}
\sum_{j=0}^p (-1)^j c_{p-j}(\omega) \, \ge\,  
\sum_{j=0}^p (-1)^j  \frac{\dim_\C H_{p-j}(M;L\otimes E)}{\dim E},
\end{eqnarray}
for $p=0, 1, 2, \dots.$ }

\begin{proof} The proof is similar to proof of Theorem 3.3. 
Let $\rho_E: \Z\pi\to \End(E_0)$ be the
monodromy representation of $E$, where $E_0$ denotes the fiber of $E$ over the base point. 
For a flat line bundle $L\to M$, the representation 
$\rho=\ml\otimes \rho_E:\Z\pi\to \End(E_0)$ defines on $E_0$ a structure of $(\C,\Z\pi)$-bimodule.

Assuming that the flat line bundle $L$
is not a $\xi$-algebraic integer (cf. above) the bimodule $E_0$ 
satisfies the condition of Theorem 1.13. 
Indeed, let $\L_0\subset E_0$ denote an integral lattice which is preserved under the monodromy
transformation along loops in $M$. Consider $Y=\Z[H]\otimes_\Z \L_0$ as a right $\Z\pi$-module,
where each $g\in \pi$ acts as follows 
\[
(h\otimes l)\cdot g =  h\cdot {\operatorname {Ab}}(g)\otimes {\operatorname {Mon}}(g)(l),\quad
  h\in H, \, l\in \L_0,\, g\in \pi,\]
where ${\operatorname Ab}:\pi\to H$ is the abelinization homomorphism and 
${\operatorname {Mon}}(g)$ is the monodromy along $g$, viewed as a map $\L_0\to \L_0$. Note that ${\operatorname {Mon}}(gg') = {\operatorname {Mon}}(g')\circ {\operatorname {Mon}}(g)$, i.e.
we indeed have a right action. $Y$ has also the obvious left $\Z[H]$ action.

Any $m\times m$ matrix $A$ with entries in $\Z\pi$
yields a $\Z[H]$-linear map $Y^m\to Y^m$ (by acting on the right), 
which may be represented by a square matrix $B$ of size
$m\cdot \dim E\times m\cdot \dim E$ with entries in $\Z[H]$. If $A$ is $\xi$-negative then $B$ is 
also $\xi$-negative.
Hence the determinant $d\in \Z[H]$
of matrix $I+B$ is an integer polynomial with $\xi$-top coefficient 1. 
Therefore, the determinant $d_L$
of the $\C$-linear map $E_0^m \to E_0^m$, given by action of $I+A$, is a nonzero complex number, 
since this determinant $d_L$ equals $\ml(d)\in \C$ and hence is nonzero, 
since otherwise we would have 
$d\in \mathcal I_L$, contradicting the assumption that $L$ is not a $\xi$-algebraic integer.

Using Theorem 1.13 we obtain that for any closed 1-form $\omega$ on $M$ having 
Morse zeros and lying in the class $\xi$, there exists a chain complex $C_\ast$ over $\C$ 
with $H_j(C_\ast) \simeq H_j(M;L\otimes E)$ and with $\dim C_j = \dim E\cdot c_j(\omega)$.
Inequalities (4.1) and (4.2) now follow from the standard argument \cite{M}.

In case when $L^\ast$ is not a $\xi$-algebraic integer, we apply the previous arguments to the class
$-\xi$ and the flat bundle $L^\ast\otimes E^\ast\otimes\mathfrak o_M$ instead of $E$ and 
use the Poincar\'e duality
as in the proof of Theorem 3.3. 
\end{proof}

\section{\bf  Generic flat vector bundles}

\subsection{} 
Let $\kk$ be a fixed algebraically closed field, possibly having a positive characteristic. 
We will consider flat $\kk$-vector bundles $E$ over a manifold $M$. 
We will understand such bundles as
locally trivial sheaves of $\kk$-vector spaces. The cohomology $H^q(M;E)$ will be understood as the
sheaf cohomology. 

Given a real cohomology class $\xi\in H^1(M;\R)$, let $\V_\xi$ be
the variety of all $\kk$-line bundles over $X$, 
which have trivial monodromy along the curves in $\ker(\xi)\subset H_1(M)$.
The variety 
$\V_\xi$ can be identified with $(\kk^\ast)^r$, where $r$ denotes the rank 
$r \, =\, \rank (H_1(M)/\ker (\xi)).$

\subsection{Definition (\cite{F3})} {\it
A flat bundle $E$ will be called $\xi$-generic if there is no $L\in \V_\xi$, so that for some $q$,
$\dim H^q(M;L\otimes E) < \dim H^q(M;E)$.}

Note that this property depends only on $\ker(\xi)$.

\subsection{Theorem} {\it Let $M$ be a closed smooth manifold and let $\xi\in H^1(M;\R)$ 
be a real cohomology class. Let $E\to X$ be a $\xi$-generic flat $\kk$-vector bundle.
Then for any closed 1-form $\omega$ on $M$ having 
Morse type zeros and lying in the class $\xi$, the number $c_p(\omega)$ of zeros
of $\omega$ having index $p$ satisfies
\[
c_p(\omega)\, \ge\,  \frac{\dim_\C H_p(M; E)}{\dim E},\qquad p=0,1,2, \dots.\]
Moreover, 
\begin{eqnarray}\,\,
\sum_{j=0}^p (-1)^j c_{p-j}(\omega) \, \ge\,  
\sum_{j=0}^p (-1)^j  \frac{\dim_\C H_{p-j}(M;  E)}{\dim E}
\end{eqnarray} for $p=0,1,2,\dots$

\begin{proof} We may view $\xi$ as a homomorphism $H_1(M;\Z)\to \R$ and let $\ker(\xi)$ denote
the kernel. Let $H$ denote $H_1(M;\Z)/\ker(\xi)$. 

We have the monodromy representation $\me: \Z\pi \to \End_{\kk}(E_0)$ of the flat bundle $E$, 
where $E_0$ is the fiber over the base point. Consider the field $\kk(H)$ of rational functions on $H$
and the following right action of $\Z\pi$ on
$X=\kk(H)\otimes_{\kk} E_0$:
\[(h\otimes l)g \, =\, h\phi(g) \otimes \me(g)(l),\quad\text{where}\quad 
g\in \pi,\, h\in \kk(H),\, l\in E_0,\]
and $\phi: \pi\to H$ is the natural projection. $X$ also has an obvious $\kk(H)$-action from the left.
Note that $X=\kk(H)\otimes_{\kk} E_0$ yields a $(\kk(H), \Z\pi)$-bimodule
satisfying the condition of Theorem 1.13 (because of the argument similar to 1.9).

Hence Theorem 1.13 applied to bimodule $X$ implies
\[
\sum_{j=0}^p (-1)^j c_{p-j}(\omega) \, \ge\,  
\sum_{j=0}^p (-1)^j  \frac{\dim_{\kk(H)} H_{p-j}(M;  X)}{\dim E}.
\]
The last inequality gives (5.1), since the assumption, that the flat bundle $E$ is 
$\xi$-generic is in fact equivalents to
\begin{eqnarray}
\dim_{\kk(H)} H_q(M;X) \, =\, \dim_{\kk}H_q(M;E), \quad\text{for any}\quad q.
\end{eqnarray}
Indeed, consider the $(\kk[H], \Z\pi)$-bimodule $Y=\kk[H]\otimes_{\kk} E_0$, defined similarly to 
$X$ but replacing the field of rational functions $\kk(H)$ by the ring $\kk[H]$ of Laurent polynomilas.
Let $C_\ast(\tilde M)$ denote the chain complex of the universal covering of $M$. For any flat 
line bundle $L\in \V_\xi$ we have the monodromy 
homomorphism $\ml: \kk[H]\to \kk$. We will view the dimension
of the homology
\[H_q(\kk\otimes_{\ml}(Y\otimes_{\Z\pi} C_\ast(\tilde M)))\simeq H_q(M;L\otimes E)\]
as a function of $L\in \V_\xi$. It is well know from 
the algebraic geometry that there exists a proper algebraic
subset $Z_q\subset \V_\xi$, so that for $L\notin Z_q$ (i.e. when $L$ is {\it a generic point})
\[\dim_{\kk} H_q(M;L\otimes E) = \dim_{\kk(H)} H_q(M;X)\]
and for $L\in Z_q$
\[\dim_{\kk} H_q(M;L\otimes E) > \dim_{\kk(H)} H_q(M;X).\]
Hence, if $E$ is $\xi$-generic, then $1\notin Z_q$ for all $q$ and therefore (5.2) follows.
\end{proof}

\section{\bf  Morse theory on manifolds with corners}

\rm 

The proof of the Main Theorem uses Morse theory for manifolds with corners. 
For this purpose we include a brief review of this theory.
The results of this section are certainly well-known, 
although it is difficult to find an appropriate reference.

\subsection{Critical points of  functions on manifolds with corners} 
Let $M$ be a $C^\infty$-smooth $n$-dimensional {\it manifold with corners}, i.e. a manifold 
which is locally diffeomorphic to $\R_+^n=(\R_+)^n$, where $\R_+$
denotes the half line $x\ge 0$. $M$ admits a natural stratification 
$M\, =\, S_0\supset S_1\supset S_2\supset\dots\supset  S_n.$ 
Here $S_1=\partial M$ is the boundary of $M$ and 
$S_k-S_{k+1}$ consists of all points $p\in M$, which have a neighborhood $U\subset M$
such that $(U,p)$ is diffeomorphic to $(\R_+^k\times \R^{n-k},0)$. We will say that the {\it points 
$p\in S_k$ are corners of order $\ge k$.}

Given a point $p\in M$, the tangent space to $M$ at $p$ is denoted 
$T_p(M)$; it is a vector space of dimension $n=\dim M$. 
We will denote by $C_p(M)$ {\it the cone of tangent directions} to $M$ at $p\in M$. By the definition,
a vector $X\in T_p(M)$ belongs to $C_p(M)$ iff there exists a smooth curve 
$\gamma: [0,\epsilon)\to M$ with $\gamma(0) =p$ and $\gamma'(0) = X$. It is clear that $C_p(M)\subset T_p(M)$ is a closed convex
cone with the property:  $C_p(M)$ is the closure of its interior in $T_p(M)$ (since any curve can be approximated by a curve with $\gamma(t)$ lying in the interior of $M$ for $t\in (0, \epsilon)$).

\subsection{Definition}{\it
Let $f: M\to \R$ be a smooth function on $M$. A point $p\in M$ is called {\it a critical point of}
$f$ if the cone of tangent directions $C_p(M)$ is disjoint from the half space
$\{X\in T_p(M); df_p(X) < 0\}$.} 

In other words, we require that at a critical point the directional derivatives
are non-negative $X(f)\ge 0$ for all vectors $X\in C_p(M)$ in the cone of tangent directions.

Let $T_p^S(M)$ denote the tangent space to the smallest stratum of $M$ containing the point $p$. 
If $p$ is a critical point
of $f$ then the differential $df_p$ vanishes on $T_p^S(X)$. The converse is not true: a point
$p\in M$ with $df_p|_{T_p^S(M)} =0$ may not be a critical point.

\subsection{Definition} {\it A critical point $p\in M$ of a smooth 
function $f:M\to \R$ is called nondegenerate (or Morse) if 

(1) $p$ is nondegenerate as a critical point of the restriction $f|_S$, 
where $S\subset M$ denotes the stratum of $M$ containing the point $p$;

(2) $C_p(M)\cap \ker[df_p: T_p(M)\to \R] = T^S_p(M)$.

Index of a nondegenerate
critical point $p\in M$ of smooth function $f:M\to \R$ 
is defined as the index of $p$ viewed as the critical point of $f|_S$.}

{\bf Examples.} 1. An interior point $p\in X$ is critical if and only if the differential $df_p$ vanishes, 
i.e. in this case the definition 1.2 coincides with the usual definition. A point $p$ lying on the boundary
$\partial X$ is a critical point of $f$ iff it is a critical point of $f|_{\partial X}$ and the derivative
$X(f)\ge 0$ is non-negative for any interior pointing tangent vector $X$.

2. Any point $p\in M$, which is a {\it local minimum point} of $f$, is a critical point of $f$.
Indeed, in the case of a local minimum point $p$ we have $C_p(M)\cap \Pi_p(f)= \emptyset$ for if
this set is nonempty we would have a curve $\gamma: [0,\epsilon)\to M$ with $\gamma(0) =p$ and $\frac{d}{dt}f(\gamma(t))<0$, which gives a contradiction.

3. A local maximum point need not to be a critical point. In fact, a local maximum
point $p\in M$ is critical if and only if $df_p=0$.

On the following
picture we see three arcs with the Morse functions given by the height ($y$-coordinate).
Boundary points A, B, D, F are Morse and have index zero. C is not a critical point. E is a degenerate
critical point.

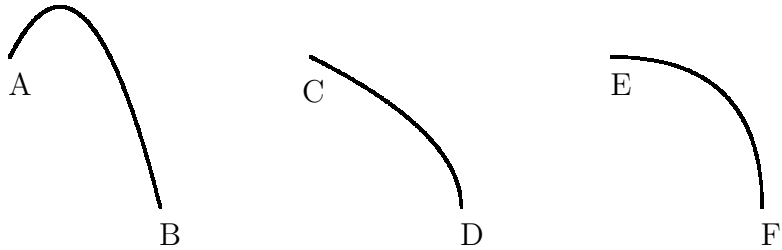
\setlength{\unitlength}{1cm}
\begin{figure}[h]
  \begin{center}
\begin{picture}(9.5,4)
\linethickness{0.3mm}
\qbezier(0,3)(1, 5)(2,1)
\qbezier(4,3)(6,2)(6,1)
\qbezier(8,3)(10,3)(10,1)
\put(0,2.5){A}
\put(2,0.5){B}
\put(3.9,2.4){C}
\put(6,0.5){D}
\put(8,2.5){E}
\put(10,0.5){F}
\end{picture}
\end{center}
\end{figure}

The following two theorems justify the above definitions. 

\subsection{Theorem} {\it Let $M$ be a manifold with corners and 
let $f: M\to \R$ be a smooth function. Suppose that the set $f^{-1}([a,b])\subset M$ is compact and contains no critical
points of $f$. Then the manifolds $f^{-1}((-\infty, a])$ and $f^{-1}((-\infty, b])$ are diffeomorphic.}

\subsection{Theorem} {\it Let $M$ be a manifold with corners and 
let $f: M\to \R$ be a smooth function. 
Suppose that the set $f^{-1}([a,b])\subset M$ is compact and 
contains a single critical point $p$, which is non-degenerate and has 
index $\lambda$.
Then the manifold $f^{-1}((-\infty, b])$ is homotopy equivalent to 
$f^{-1}((-\infty, a])\cup e^\lambda$,
the result of glueing a cell of dimension $\lambda$ to $f^{-1}((-\infty, a])$.}

The proofs are obtained by
the usual arguments 
of the Morse theory using the gradient flows.

\section{\bf  Chain collapse}

In this subsection we remind a technique, developed in \cite{FR}, which is a chain
analog of the well-known geometric operation of {\it combinatorial collapse}. The later
is illustrated by the
following picture: a {\it "protruding"}
 cell $e'$, which is not a part of the boundary of any other cell, has a
{\it "free"} face $e$, and we may remove the cells $e'$ and $e$ without changing the 
(simple) homotopy type.

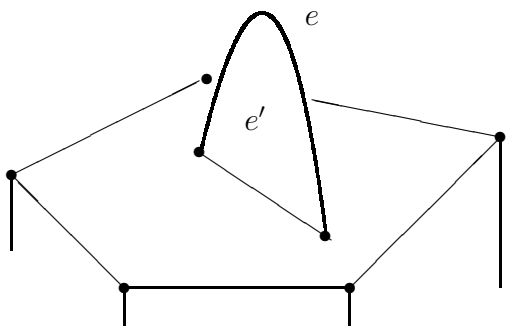
\begin{figure}[h]
  \begin{center}
\begin{picture}(9,5)
\linethickness{0.2mm}
\put(3,1){\line(1,0){3}}
\put(3,1){\line(-1,1){1.5}}
\put(6,1){\line(1,1){2}}
\put(1.5,2.5){\line(2,1){2.5}}
\put(8,3){\line(-5,1){2.5}}
\put(3,1){\circle*{0.15}}
\put(1.5,2.5){\circle*{0.15}}
\put(8,3){\circle*{0.15}}
\put(4,2.8){\circle*{0.15}}
\put(6,1){\circle*{0.15}}
\put(5.67,1.69){\circle*{0.15}}
\put(4.1,3.785){\circle*{0.15}}
\put(3,1){\line(0,-1){0.5}}
\put(6,1){\line(0,-1){0.5}}
\put(1.5,2.5){\line(0,-1){1}}
\put(8,3){\line(0,-1){2}}
\put(4.6,3.1){$e'$}
\put(5.4,4.5){$e$}
\put(4.1,3.785){\line(5,-1){0.3}}
\linethickness{0.3mm}
\put(4,2.8){\line(3,-2){1.75}}
\qbezier(4,2.8)(5,7)(5.67,1.69)
\end{picture}
\end{center}
\end{figure}

In the sequel we will use the following matrix convention.
A morphism between direct sums of modules
\[\phi~=~(\phi_{ij})~:~P_1 \oplus P_2 \oplus \dots \oplus P_n
\to Q_1 \oplus Q_2 \oplus \dots \oplus Q_m\]
is denoted by an $m \times n$ matrix with entries morphisms
$\phi_{ij}:P_j \to Q_i$, so that
\[\phi(x_1,x_2,\dots,x_n)~=~(\sum\limits^n_{j=1}\phi_{1j}(x_j),
\sum\limits^n_{j=1}\phi_{2j}(x_j),\dots, \sum\limits^n_{j=1}\phi_{mj}(x_j))~.\]

The following Lemma is a minor generalization of Lemma 2.3 of \cite{FR}.

\subsection{Lemma (Chain analog of combinatorial collapse)}{\it
Let $\RR$ be a ring and let $(B,d_B)$ be a chain complex of free $\RR$-modules 
with 
$$B_i~=~D'_i \oplus D_i\oplus C_{i}, \quad i=0, 1, 2, \dots$$ 
and with the differential $d_B$ having the form
\begin{eqnarray}
d_B~=~\left[\matrix d_{D'} & 0 & 0 \cr 
\gamma & d_D & \alpha \cr
\beta & \sigma & d_{C} \endmatrix\right]\end{eqnarray}
where $\alpha: C_i\to D_{i-1}$, $\beta: D'_i\to C_{i-1}$, $\gamma: D'_i\to D_{i-1}$ and $\sigma: D_i\to C_{i-1}$.
Suppose that the homomorphism $\gamma: D'_i\to D_{i-1}$ is an isomorphism for all $i$.
Then formula
\begin{eqnarray}
\hat{d}_C~=~d_C-\beta \gamma^{-1}\alpha~:~C_i \to C_{i-1}
\end{eqnarray}
defines a "deformed differential" on $C$ (i.e. ${\hat{d}_{{}_C}{}}^2=0$), and the chain complexes $(B,d_B)$ and 
$(C, \hat{d}_C)$ are chain homotopy equivalent.}

The Lemma states that one may perform a deformation (given by (7.2)) to achieve cancelation
of the chain complex. The role of the "protruding" 
cell $e'$ plays the submodule $D'$. It does not appear as a part of boundary of any other chain (we see two zeros in the first row of the matrix presenting the 
differential $d_B$). The role of the free face $e$ plays $D$; the fact that $\gamma: D'_i\to D_{i-1}$
is an isomorphism expresses that.

\begin{proof} First note that ${(d_B)}^2=0$ implies
$${(d_D)}^2+\alpha\sigma =0,\quad {(d_C)}^2+\sigma\alpha=0,\quad {(d_{D'})}^2=0,$$
and also

\begin{eqnarray}
d_D\alpha +\alpha d_C&=&0,\nonumber\\
d_C\beta +\beta d_{D'}+\sigma\gamma &=&0,\\
d_D\gamma +\gamma d_{D'}+\alpha\beta \, &=&\, 0,\nonumber\\
d_C \sigma +\sigma d_D &=&0.\nonumber
\end{eqnarray}
Using (7.3) we obtain
$$
\aligned
{\hat{d}_{{}_C}{}}^2&=(d_C -\beta\gamma^{-1}\alpha)\cdot (d_C -\beta\gamma^{-1}\alpha) =\\
&= -\sigma\alpha - d_C\beta\gamma^{-1}\alpha - \beta\gamma^{-1}\alpha d_C + 
\beta\gamma^{-1}\alpha\cdot\beta\gamma^{-1}\alpha=\\
&=-\sigma\alpha +[\sigma\gamma +\beta d_{D'}]\gamma^{-1}\alpha +
\beta\gamma^{-1}d_D\alpha - 
\beta\gamma^{-1}[d_D\gamma +\gamma d_{D'}]\gamma^{-1}\alpha=0.
\endaligned
$$

Now we define two chain maps:
\[f\, =\, \left[\matrix 0 & -\beta\gamma^{-1} & 1   \endmatrix\right]: B \to C,\quad
\text{and}\quad g\, =\, \left[\matrix -\gamma^{-1}\alpha\\
 0\\
1\endmatrix\right]: C\to B.\]
One checks that
\[fd_B = \hat{d_C}f,\quad d_B g = g \hat{d_C}, \quad fg = 1_C\]
and 
$$gf = 1_B - d_B h -hd_B,\quad\text{where}\quad h\, =\, \left[\matrix 0 & \gamma^{-1} & 0\\
0 & 0 & 0\\
0 & 0 & 0\endmatrix\right]: B\to B.
$$
Hence $f$ and $g$ are mutually inverse homotopy equivalences.

This completes the proof. 
\end{proof}

\subsection{} Under the conditions of Lemma 7.1
we will say that the chain complex $(B, d_B)$ {\it collapses} to the complex
$(C, \hat{d_C})$.

For our purposes in this paper 
the most important will be the special case of Lemma 7.1, when $\alpha=0$.
The formula (7.2) for the differential $\hat{d_C}$ shows that in that case we will have $\hat{d_C}=d_C$, i.e. the differential of the collapsed complex is given by the original matrix, in 
which we simply ignore the terms contained in $D'\oplus D$. 
In this case we will say that the chain complex
$(B, d_B)$ {\it simply collapses} to the complex
$(C, d_C)$.

Simple collapse is a full algebraic analog of the
elementary combinatorial collapse shown on the picture in the beginning of \S 7.

\section{\bf  Proof of Theorem 1.11}

The proof will consist in the following steps. Firstly, we describe a cell decomposition 
of the manifold $M$ (up to homotopy type)
related to the given closed 1-form $\omega$. 
Secondly, we compute the chain complex of the
universal covering, corresponding to this cell decomposition. Thirdly, we apply the chain collapse
technique, allowing cancelation of the chain complex after a suitable Cohn localization, which was described in section \S 7.

\subsection{} Let $\omega$ be a given closed 1-form on an $n$-dimensional manifold $M$ 
having only Morse zeros. Represent $\omega$ as a linear combination
\begin{eqnarray}
\omega \, =\, \sum_{i=1}^r \lambda_i\omega_i,\quad \lambda_i >0, \quad \lambda_i\in \R,
\end{eqnarray}
of closed 1-forms $\omega_i$, 
so that the cohomology classes $[\omega_1], \dots, [\omega_r]\in H^1(M;\R)$ are integral 
(i.e. $[\omega_i]\in H^1(M;\Z)$), indivisible, and linearly independent. 

\subsection{Lemma} {\it Representation (8.1) can be chosen so 
that for any $p\in M$, which is not a zero of the form $\omega$, the intersection of the half-spaces
\[
\{X\in T_p(M); \omega(X)>0\} \cap\bigcap_{j=1}^r\{X\in T_p(M); \omega_j(X)>0\}\]
is nonempty.}

\begin{proof} Fix on $M$ a Riemannian metric and let $X$ be the gradient field of $\omega$ 
with respect to this metric. 
If $U\subset M$ denotes the union of open $\delta$-balls ($\delta>0$ is small)
around the zeros $p_1, \dots, p_k$ of $\omega$, then on $M-U$ holds
$\omega(X)>\epsilon>0$ for some $\epsilon >0$. Hence for any closed 1-form $\tilde\omega$,
which sufficiently closely approximates $r^{-1}\omega$ on the compact subset $M-U$, holds
$\tilde \omega(X)>0$ on $M-U$. It is clear that we may find representation (8.1) with the forms 
$\lambda_i\omega_i$ approximating $r^{-1}\omega$ on $M-U$ and such that 
on each of the $\delta$-balls the forms $\lambda_i\omega_i$ coincide with $r^{-1}\omega$. 
Then the gradient $X$ of $\omega$ will belong to the intersection of all the half-spaces above
for all $p\in M$, distinct from $p_1, \dots, p_k$.

\end{proof}
\subsection{} For any $i=1, 2, \dots, r$ there exists a smooth map
$f_i: M\to S^1$, so that $\omega_i=f_i^\ast(d\theta)$, where $d\theta$ is the angular form on the 
circle $S^1$. We obtain a smooth map from $M$ to an $r$-dimensional torus
$f=(f_1, f_2, \dots, f_r):  M \to T^r = S^1\times S^1\times \dots \times S^1.$
Choose a {\it generic} point $a_i\in S^1$ for each $i=1, \dots, r$.
Denote $V_i=f_i^{-1}(a_i)$, where $i=1, 2, \dots, r$. These are smooth codimension 1 
submanifolds of $M$, which meet transversally. Each $V_i$ naturally comes 
equipped with a specified 
orientation of its normal bundle: a vector $X\in T_p(M)$, where $p\in V_i$, is {\it positive}
if $\omega_i(X)>0$.

Let manifold $N$ be obtained from the manifold $M$ by cutting
along the submanifolds $V_1, V_2, \dots, V_r$. $N$ is a manifold with corners.
We will denote the parts of its boundary as follows. The manifold 
\[V_i - \bigcup_{j\ne i} V_j\]
determines two diffeomorphic $(n-1)$-dimensional submanifolds of $\partial N$. 
Their closures we will denote by $N_i^+, N_i^-\subset \partial N$. They are 
canonically diffeomorphic
$J_i: N_i^+\to N_i^-$; the points $x\in N_i^+$ and $J_i(x) \in N_i^-$ are 
mapped to the same point when
we glue $M$ back from $N$. The notations $N_i^\pm$ are chosen so that the positive 
(with respect to
$\omega_i$) normal field points inside $N$ along $N_i^+$ and points outside $N$ along $N_i^-$.

Let $\Pi: N\to M$ denote the canonical identification map. The induced 
form $\Pi^\ast\omega$ is exact. Also, the forms $\Pi^\ast\omega_i$ are exact, $i=1, \dots, r$.
Hence we obtain smooth functions $g_i: N\to \R$, $i=1, 2, \dots, r$, so that
$f_i(\Pi(x)) = \exp(2\pi i g_i(x)), \quad x\in N.$

\subsection{} We obtain the following Morse function 
\[
g\, = \, \sum_{i=1}^r \lambda_i g_i \, : \, N \, \to \, \R, \quad dg = \Pi^\ast\omega.
\]
Consider the critical points of $g$, viewed 
as a Morse function on manifold with corners $N$, cf. \S 6.
It is clear that if $x\in N$ is a critical point of $g$, belonging to the interior $\In (N)$ of $N$, 
then the point $\Pi(x)\in M$ is a zero of the form $\omega$, and, conversely, zeros of $\omega$ produce the internal critical points of $g$. Also, the index 
of $x\in \In(N)$ as a critical point of $g$
is the same as the index of $\Pi(x)$ as a zero of $\omega$. 

Let us show that the function $g$
has no critical points on the "negative part of the boundary" 
\[\partial_- N \, =\, \partial N - \bigcup_{i=1}^r N_i^+.\]
Indeed, if a point $p\in \partial_- N$ belongs to intersection 
$N_{i_1}^-\cap N_{i_2}^-\cap\dots \cap
N_{i_k}^-$ and is a corner of order $k$, then the cone
of tangent directions $C_p(N)$ (cf. \S 6) coincides with the intersection of the half-spaces
\[\bigcap_{s=1}^k \{X\in T_p(N); X(g_{i_s})< 0\}.\]
Hence, by Lemma 8.2, we may find $X\in C_p(N)$, so that $X(g)<0$, and thus $p$ is not a critical 
point of $g$, according to the definition given in \S 6. 

As the result, we obtain (using Theorem 6.2) that homotopy type of $N$ is obtained by glueing cells corresponding to the boundary critical points of $g$, lying in the positive part of the boundary
$\cup N_j^+$, and also to the internal critical points of $g$; the later are in 1-1 
correspondence with the zeros of $\omega$. 

\subsection{Cell decomposition of $N$} 
To simplify the notations we will denote the submanifold $N_j^+$ by $N_j$. Any intersection
$N_{i_1}\cap N_{i_2}\cap\dots\cap N_{i_k}$, for a sequence
$i_1< i_2< \dots i_k\le r$ of indices, will be denoted by $N_{i_1 i_2\dots i_k}$ or even
shorter by $N_\alpha$, where $\alpha =\{i_1, i_2, \dots, i_k\}$ is {\it a multi-index}, i.e.
a subset $\alpha\subset \{1, 2, \dots, r\}$. 
The number $k=|\alpha|$ will be called {\it the length} of the multi-index $\alpha$. 
The submanifold $N_\alpha$ consists of points where $N$ has
corner of order $\ge k$, cf. \S 6. It will be convenient to allow also the empty subset
$\alpha=\emptyset$ as a milti-index; we will understand $N_{\emptyset}$ as $N$.

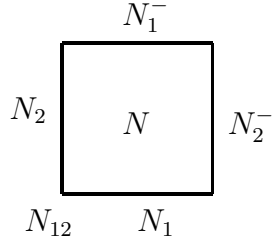
\begin{figure}[h]
  \begin{center}
\begin{picture}(9,4)
\linethickness{0.3mm}
\put(3,1){\line(1,0){2}}
\put(3,1){\line(0,1){2}}
\put(3,3){\line(1,0){2}}
\put(5,1){\line(0,1){2}}
\put(2.5,0.5){$N_{12}$}
\put(4,0.5){$N_1$}
\put(3.8,1.8){$N$}
\put(2.3,2){$N_2$}
\put(3.8,3.25){$N_1^-$}
\put(5.2,1.8){$N_2^-$}
\end{picture}
\end{center}
\end{figure}

Suppose that a multi-index $\beta$ is obtained
by removing $i$ from a multi-index $\alpha$.Then there is an inclusion
$L_i: N_\alpha \to N_\beta$. For example, we have an obvious inclusion $L_j$ of $N_j$ into 
$N=N_{\emptyset}$ and also an inclusion
of $N_{ij}=N_i\cap N_j$ into $N_i$. We have the commutativity relations
\[L_i\circ L_j = L_j\circ L_i, \quad i\ne j.\]

Fix a cell decomposition of $\cup_{i=1}^r N_i$ with the following property: $N_\alpha$ is 
a subcomplex of $N_\beta$ for any pair of multi-indices $\alpha =(i_1, i_2, \dots, i_k)$ and 
$\beta= (j_1, j_2, \dots,j_l)$, assuming that $l<k$ and all indices $j_s$ appear also in $\alpha$.
It clear that such cell decomposition exists. One simply starts with the corners of higher order and
acts inductively. 

A cell $e$ of $N_\alpha$ will be called {\it proper} 
if it is disjoint from all $N_\beta\subset N_\alpha$ with $|\beta|>|\alpha|$.

A cell decomposition of $N$, up to homotopy, is obtained by adding to the cell decomposition of 
$\cup_{i=1}^r N_i$,
described above, the cells corresponding to the zeros of the form $\omega$, i.e. the 
internal critical points of $g$. 

\subsection{Cell decomposition of $M$} Manifold $M$ is obtained from $N$ by 
identification of
faces according to the diffeomorphisms
$J_i: N_i \to N_i^-\subset N$, where $i=1, 2, \dots, r$. 
In fact, $J_i$ acts also on manifolds $N_\alpha$: 
if milti-index $\alpha$ contains some index $i\in \{1, 2, \dots, r\}$, then diffeomorphism 
$J_i$ maps the manifold 
$N_\alpha$ into $N_\beta$, where $\beta$ is obtained from $\alpha$ by removing $i$. 
The maps $J_i$ commute with each other
\[J_i\circ J_j = J_j\circ J_i, \quad i\ne j.\]
Also, we obviously have
\[L_i\circ J_j = J_j\circ L_i,\quad i\ne j,\]
i.e. the maps $J_j$ commute with the inclusions $L_i$.

For a subset $\alpha\subset \{1, 2, \dots, r\}$ we will denote by $I^\alpha$ the set of functions
from $\alpha$ to the interval $I=[-1,1]$ (cube of dimension $k=|\alpha|$). 
Given an index $i\in \alpha$,
we will consider the following face inclusions
$\Psi^\pm_i: I^{\beta}\, \to \, I^{\alpha},$ where $\beta=\alpha-\{i\},$
which act as follows: for a function $t: \beta\to I$, its image $\Psi^\pm_i(t)$ is defined
by  
\[\Psi_i^\pm|_\beta = t, \quad \Psi_i^\pm(i) = \pm 1.\]

\subsection{Lemma}{\it
The manifold $M$ is homeomorphic to the factor-space of the disjoint union
\begin{eqnarray}\bigcup_{\alpha} N_\alpha\times I^{\alpha}\, /\sim\, ,\end{eqnarray}
where $\alpha$ runs over all multi-indices (including the empty set $\alpha = \emptyset$), 
and where $\sim$ is the following equivalence relation.
For any point $a\in N_\alpha$ with $|\alpha|=k$, and for any $t\in I^{\beta}$, where 
$\beta =\alpha -\{i\}$, we identify the following
pairs of points of the union (8.2)}
\begin{eqnarray}(a, \Psi_i^+(t)) \, \sim\,  (L_i(a), t),\quad i\in \alpha,\end{eqnarray}
\begin{eqnarray}(a, \Psi_i^-(t)) \, \sim \, (J_i(a), t),\quad i\in \alpha.\end{eqnarray}

For example,  a simple face $N_i$ contributes a cylinder $N_i\times I$ and for $a\in N_i$ we 
identify $(a,1)$ with $L_i(a)=a \in N$ and glue $(a, -1)$ to $J_i(a)$. 
Hence two ends of the cylinder $N_i\times I$ are glued to $N$.

\subsection*{Proof of Lemma} For any $i=1, 2, \dots, r$, consider a parallel copy $V_i'$ of the 
submanifold
$V_i$ shifted slightly in the positive normal direction. 
Intersections of the submanifolds $V'_i$
with $N$ produce parallel copies of the faces $N_i$, shifted slightly up (shown by dotted lines
on the left picture). If we cut $M$ along the submanifolds $V_1$, $V_2, \dots, V_r$ and 
also along $V'_1, V'_2,\dots, V'_r$ we obtain what is shown on the right picture. 
The large rectangle represents
a manifolds with corners, which is diffeomorphic to $N$. Each of the small rectangles is 
diffeomorphic to $N_\alpha\times I^{\alpha}$ for a suitable multi-index $\alpha$. 
The small rectangles are glued to $N$ and to
each other in accordance with relations (8.3), (8.4). 

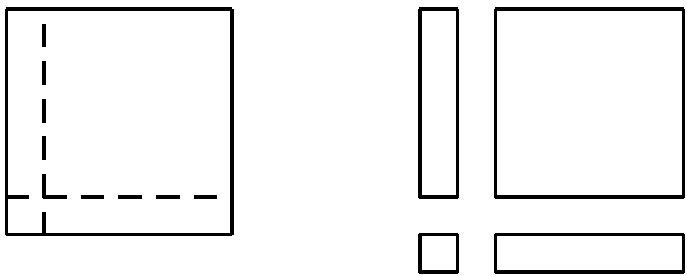
\begin{figure}[ht]
  \begin{center}
\begin{picture}(13,4.5)
\linethickness{0.3mm}
\put(1,1){\line(1,0){3}}
\put(1,1){\line(0,1){3}}
\put(4,1){\line(0,1){3}}
\put(1,4){\line(1,0){3}}
\multiput(1,1.5)(0.5,0){6}{\line(1,0){0.3}}
\multiput(1.5,1)(0,0.5){6}{\line(0,1){0.3}}

\put(7.5,1.5){\line(1,0){2.5}}
\put(7.5,1.5){\line(0,1){2.5}}
\put(10,1.5){\line(0,1){2.5}}
\put(7.5,4){\line(1,0){2.5}}

\put(6.5,1.5){\line(1,0){0.5}}
\put(6.5,1.5){\line(0,1){2.5}}
\put(7,1.5){\line(0,1){2.5}}
\put(6.5,4){\line(1,0){0.5}}

\put(7.5,1){\line(1,0){2.5}}
\put(7.5,0.5){\line(0,1){0.5}}
\put(10,0.5){\line(0,1){0.5}}
\put(7.5,0.5){\line(1,0){2.5}}

\put(6.5,0.5){\line(1,0){0.5}}
\put(6.5,0.5){\line(0,1){0.5}}
\put(7,0.5){\line(0,1){0.5}}
\put(6.5,1){\line(1,0){0.5}}

\end{picture}
\end{center}
\end{figure}

This shows, that the factor-space of the disjoint union 
$\bigcup_{\alpha} N_\alpha\times I^{\alpha}$ 
with respect to relations (8.3) alone, gives a manifold homeomorphic to $N$. 
Now, identifying points of the obtained space
according to relations (8.4) precisely means glueing each face $N_i$ of $N$ 
to its image under $J_i$.

Now we may describe the cell structure of $N$.
Each $N_\alpha\times I^{\alpha}$, where $\alpha\ne \emptyset$, 
has an obvious cell structure, since $N_\alpha$ is a subcomplex of $\cup_{i=1}^r N_i$ and 
$I^{\alpha}$ has the product cell structure (the interval $I$ is assumed to be decomposed into
union of two its end points and the open interval). The subcomplex
$N_\alpha\times \partial(I^{\alpha})$ is identified via relations (8.3), (8.4) with the 
points determined by the  
lower codimension strata. Hence each product $N_\alpha\times I^{\alpha}$, where $\alpha\ne\emptyset$,
really adds cells
of the form $e\times \In (I)^{k}$ for any cell $e\subset N_\alpha$, where $k=|\alpha|$. 
We will denote the cell $e\times \In (I)^{k}$ by symbol $Q_\alpha(e)$. 

Note that any cell $e\subset N_\alpha$ is contained also in all $N_\beta$, where $\beta$ is a subset
of $\alpha$. The resulting cells $Q_\alpha(e)$ and $Q_\beta(e)$ will be, of course, distinct.

Hence, we obtain that a cell complex $Y$, homotopy equivalent to $M$, can be obtained as follows.
Start with a cell decomposition of $\cup_{i=1}^r N_i$ as above. Glue to it the cells, 
corresponding to the internal critical points of $g$. The result will give a cell complex
$N'$ containing the union $\cup_{i=1}^r N_i$ and homotopy equivalent to 
$N=N_{\emptyset}$ relative to $\cup_{i=1}^r N_i$.
After that for any proper cell $e\subset N_\alpha$, where
$\alpha\ne \emptyset$,
introduce 
cells $Q_\beta (e)$ (described above), where $\beta$ runs over all 
sub-indices of $\alpha$, so that the dimension of the cell
$Q_\beta (e)$ equals 
\begin{eqnarray}\dim (Q_\beta (e))\, =\, \dim e + |\beta|.\end{eqnarray}
Note that $Q_\emptyset(e)$ is just $e$. 

Hence, each proper cell $e\subset N_\alpha$ produces $2^{k}$ cells, where $k=|\alpha|$,
(including $e$ itself)
of the CW-decomposition 
$Y$. 

\subsection{The chain complex}
Here we will calculate the cellular chain complex $C_\ast(\tilde Y)$ of the universal covering $\tilde Y$
of the complex $Y$ described in \S 8.7. 

Each cell of $Y$ can be lifted into the covering $\tilde Y$ and all possible lifts are parameterized 
by the elements of the fundamental group $\pi =\pi_1(Y, v_1)= \pi_1(M)$.
It will be convenient to make a choice of a base point $v_1\in Y$,
so that $v_1\notin \cup_{i=1}^r N_i \subset Y$. 

In order to fix the lifts of the cells of $Y$ we will describe for each cell of $Y$ a path in $Y$, starting from the base point $v_1$ and leading to an internal point of the cell. We will call this path 
{\it the tail of} $e$.
After an arbitrary choice of a lift of the base point $v_1$, we will obtain lifts of all the cells into the covering $\tilde X$, determined (in an obvious way) by the choice of the tails.
In general the complex $N'$ may contain several connected components. 
We will assume now that $N'$ is connected; the case when $N'$ is disconnected will be treated
later in 8.11.

Now we will describe the choice of the tails of the cells of $Y$.
We will refer to the cell structure of $Y$ described above.
For any cell $e$ of complex $N'\subset Y$ we will
choose the tail $\gamma_e:[0,1]\to N'$ as an arbitrary path in $N'$, joining the base point 
with an internal point of $e$.
The tails of the cells of the form $Q_\beta(e)$, where $\beta\ne \emptyset$, $\beta\subset\alpha$, and 
$e\subset N_\alpha$ is a proper cell, are constructed by induction of the length $|\beta|$. Indeed, 
let $\beta=\{i_1, i_2, \dots, i_m\}$, where $i_1<i_2<\dots < i_m$. Then the cell $Q_\beta(e)$ 
by one of its faces is glued to $Q_{\beta-i_1}(e)$ (cf. (8.3)), 
and we will fix the tail of $Q_\beta(e)$, which
first travels from the base point to an internal point of $Q_{\beta-i_1}(e)$ along the tail of
$Q_{\beta-i_1}(e)$ and then drops slightly
into an internal point of $Q_\beta(e)$.

We will assume that the cells $e\subset N$ are oriented. Then each cell $Q_\beta(e)$ (being 
$e\times \In(I^\beta)$) gets an orientation. Here we assume that the orientation of the cube
$I^\beta$ is canonical, i.e. the product orientation 
$I^\beta = I^{i_1}\times I^{i_2}\times \dots \times I^{i_m}$,
where $\beta=\{i_1, i_2, \dots, i_m\}$ and $i_1<i_2<\dots <i_m$.

Consider the cellular chain complex $C_\ast(\tilde Y)$. The cells $Q_\beta(e)$, together with their
tails and the orientations, 
define a free $\Z\pi$-basis of $C_\ast(\tilde Y)$. In order to distinguish the cells from their
corresponding basis elements, we will denote the later by $P_\beta(e)$.

\subsection{Lemma} {\it Let $P_\beta(e)$ be a basis element,
where $e\subset N_\alpha$ is a proper cell and $\beta$ is a non-empty subset, 
$\beta=\{i_1, i_2, \dots, i_m\}\subset \alpha$
with $i_1<i_2<\dots < i_m$. 
The action of the boundary homomorphism of the complex $C_\ast(\tilde Y)$ 
on $P_\beta(e)$ is given by the following formula:
\begin{eqnarray}
\aligned
\quad \partial P_\beta(e)\, & = \, P_\beta (\partial e) + \\
&+\sum_{s=1}^m (-1)^{\dim e + s+1} \cdot [P_{\beta -\{i_s\}}(e) \, 
-\sum_{\stackrel {e'\subset N_{\alpha-\{i_s\}}}{ \dim e' =\dim e}}
 \langle e:e'\rangle_{i_s}\cdot P_{\beta-\{i_s\}}(e')].
\endaligned
\end{eqnarray}
In the second sum $e'$ runs over all the cells of 
$N_{\alpha-\{i_s\}}$ having dimension $\dim e$ and the coefficient 
$\langle e:e'\rangle_{i_s}\in \Z\pi$ denotes an element of the group ring, 
which is independent of $\beta$ and is
$\xi$-negative (cf. 1.3), where $\xi=[\omega]:\pi \to \R$ is the cohomology class of $\omega$.}

The notation $P_\beta (\partial e)$ appearing in (8.6), is understood as follows. 
$\partial e$ can be represented as a sum
$\sum c_j e_j$, where $e_j$ are cells of $N$ and $c_j\in \Z$. Then $P_\beta (\partial e)$ is defined
as $\sum c_j P_\beta (e_j)$. 

Note that in case $|\alpha|=1$ the second sum in (8.6) includes basis elements of the form 
$P_\emptyset(e') =e'$, where $e'$ is a cell on $N=N_\emptyset$, i.e. the cells 
originating from the critical points of the 1-form $\omega$, lying in the interior of $N$.

\subsection*{Proof of Lemma} Let $B$ be a ball of dimension $\dim e$ and let 
$\chi: B\to N_\alpha$ be the characteristic map of the cell $e$. Then the
characteristic map of the cell $Q_\beta(e)\subset Y$ is the composite
\[\chi_\beta\, : \, \, B\times I^\beta \, \, \, 
\stackrel{\chi\times 1} \longrightarrow\,  N_\alpha\times I^\beta  
\longrightarrow  N_\beta\times I^\beta\, \longrightarrow Y,\]
where $N_\alpha \to N_\beta$ is the inclusion. The boundary of the disk 
$B\times I^\beta$ is the union
$$(\partial B) \times I^\beta \, \cup \, \bigcup_{s=1}^m \, B\times \Psi^\pm_{i_s}(I^{\beta-\{i_s\}})$$
and $\chi_\beta$, mapping each part of this decomposition, contributes a corresponding term into
formula (8.6).

The map $\chi_\beta: (\partial B) \times I^\beta\to Y$ clearly represents the chain 
$P_\beta(\partial e)\in C_\ast(\tilde Y)$. 

Because of
identifications (8.3), the map
\[\chi_\beta: \, B\times \Psi^+_{i_s}(I^{\beta-\{i_s\}})\to Y\] 
represents the chain 
$(-1)^{\dim e +s +1}P_{\beta-\{i_s\}}(e)\in C_\ast(\tilde Y)$. 

Now we want to examine the chain in $\tilde Y$ determined by restricting the map $\chi_\beta$
onto 
$B\times \Psi^-_{i_s}(I^{\beta-\{i_s\}})$. Because of identification (8.4),
it can also be represented by the map
\begin{eqnarray}\quad\quad\quad
B\times I^{\beta-\{i_s\}}\, \stackrel {\chi\times 1} \longrightarrow\, 
N_\alpha\times I^{\beta-\{i_s\}}\to 
N_\beta\times I^{\beta-\{i_s\}}\stackrel{J_{i_s}\times 1}\longrightarrow
N_{\beta-\{i_s\}}\times I^{\beta-\{i_s\}}.\end{eqnarray}
Denote by
\begin{eqnarray}
\sum_{\stackrel {e'\subset N_{\alpha-\{i_s\}}} {\dim e'=\dim e}}\,  
\langle e:e'\rangle_{i_s}\cdot e', \quad\text{where}\quad 
\langle e:e'\rangle_{i_s}\, \in\, \Z\pi,\end{eqnarray}
the chain in $C_\ast(\tilde N_{\alpha-\{i_s\}})$,
determined by the map 
\[\chi_s:\, B\stackrel \chi \longrightarrow N_\alpha \stackrel {J_{i_s}} 
\longrightarrow N_{\alpha-\{i_s\}},\]
where $\tilde N_{\alpha-\{i_s\}}$ denotes the part of the universal covering $\tilde Y$ lying over $N_{\alpha-\{i_s\}}$.
We will specify the tail of cell $\chi_s(B)\subset N_{\alpha-\{i_s\}}$ as follows: it first travels
along the tail of the cell $e\subset N_\alpha$ reaching a point $x\in N_\alpha$, which is identified
via (8.4) with $J_{i_s}(x)\in N_{\alpha-\{i_s\}}$.

The chain in $C_\ast(\tilde Y)$, represented by (8.7), is obtained by applying the operator \newline
$(-1)^{\dim e+s}P_{\beta-\{i_s\}}$ to (8.8). 
The result coincides with the second sum in (8.6).

Hence, to finish the proof, we only need to show that
the incidence coefficients $\langle e:e'\rangle_{i_s}\in \Z\pi$, which appear in (8.8), are $\xi$-negative.
We may assume that the map $\chi_s: B\to N_{\alpha-\{i_s\}}$ takes values in the 
$(\dim e)$-dimensional skeleton and 
maps the interior of 
$B$ smoothly on the cells $e'\subset N_{\alpha-\{i_s\}}$. Then, the incidence coefficient
$\langle e:e'\rangle_{i_s}\in \Z\pi$ can be computed as follows. Fix a point $p'\in e'$ such that $\chi_s$
is transversal with respect to $p'$. The preimage $\chi_s^{-1}(p')$ consists of finitely many points
$\{p_1, \dots, p_l\}\subset B$. The incidence coefficient equals
\begin{eqnarray}\langle e:e'\rangle_{i_s} \, =\, \sum_{j=1}^l \sign(p_j)[\gamma_j]\, \in \Z\pi,\end{eqnarray}
where $\sign(p_j)=\pm 1$ is a sign determined by the orientations of the cells $e$ and $e'$, and
$\gamma_j$ is a loop in $Y$ which is defined as follows. First, we will extend
the tails of the cells 
$e$ and $e'$ by arcs lying inside the cells, such that the points $\chi_s(p_j)\in e$ and 
$p'\in e'$ be the end points of the corresponding tails. Then $\gamma_j$ is the loop, which
starts at the base point, follows the tail of $e$, and then returns back along the tail of $e'$, 
moving in the reverse direction. Since the
tail of $e$ makes once a jump $J_{i_s}$ (crossing the submanifold $N_{i_s}$ in the negative direction)
we obtain that the cohomology classes of the forms $\omega_1, \dots, \omega_r$ take on the loop
$\gamma_j$ the following values:
$$
\langle[\omega_i], [\gamma_j]\rangle \,\, =\, \,  \cases -1, \quad\text{for}\quad i=i_s,\\
0,\quad\text{for}\quad  i\ne i_s.\endcases
$$
Hence we obtain that
\[
\langle\xi, [\gamma_j]\rangle \, \, \, \le\, \, \, -\lambda <0, 
\quad\text{where}\quad \lambda=\min_{1\le j\le r}\lambda_j.
\]
This completes the proof of the Lemma. \qed

\subsection{} One may rewrite formula (8.6) as follows
\begin{eqnarray}
\aligned
\partial P_\beta(e)\, & = \, P_\beta (\partial e) + \\
&+\sum_{s=1}^m (-1)^{\dim e + s+1} \cdot (1-\langle e:e\rangle_{i_s}) \cdot P_{\beta -\{i_s\}}(e) +\\
&+\sum_{s=1}^m\sum_{{\stackrel {e'\ne e} {\stackrel {e'\subset N_{\alpha-\{i_s\}}} { \dim e' =\dim e}}}}
 (-1)^{\dim e +s} \cdot \langle e:e'\rangle_{i_s}\cdot P_{\beta-\{i_s\}}(e').
\endaligned\end{eqnarray}

\subsection{} Now we assume that $N$ is disconnected and show how one may 
choose the tails of the cells of $Y$ so that the statement of Lemma 8.9 remains true.

First we choose base points $w_1, w_2, \dots, w_l\in \In(N)$, one in each connected component of
$N$. Consider their images $\Pi(w_1), \dots, \Pi(w_l)\in M - \bigcup_{i=1}^r V_i$ under the identification map $\Pi: N\to M$, cf. 8.3. Let us show that {\it one may find smooth paths
$$\mu_j: [0,1]\to M, \, \, \mu_j(0) = \Pi(w_1), \, \, \mu_j(1) = \Pi(w_j), \, \, j=1, 2, \dots, l,$$
such that
\begin{eqnarray}
|\sum_{i=1}^r \lambda_i (\mu_j\cdot V_i)| \, <\, \frac{1}{3}\lambda,\quad\text{where}\quad
\lambda=\min\{\lambda_1, \dots, \lambda_r\}.\end{eqnarray}}
Here $\mu_j\cdot V_i\in \Z$ denotes the intersection number of $\mu_j$ and $V_i$ (recall that the
normal bundles to the submanifolds $V_i$ are oriented).

Indeed, in case $r=1$ we have only one submanifold $V_1\subset M$ and there exists a closed
loop $\delta$ in $M$, starting and ending at $w_1$, 
with $\delta\cdot V_1=1$ (since the cohomology class $[\omega]\in H^1(M;\Z)$
is integral and indivisible). Hence if we have an arbitrary path $\overline\mu_j$ in $M$ 
joining the points 
$\Pi(w_1)$ and $\Pi(w_j)$, then the path 
$\mu_j= \overline\mu_j - (\overline\mu_j\cdot V_1)\cdot \delta$ 
(i.e. going several times along the loop $\delta$ and then along $\overline\mu_j$) 
satisfies our requirements (8.11).

In the case $r>1$ the numbers $\lambda_1>0, \lambda_2>0, \dots, \lambda_r>0$ are linearly
independent over $\Z$ and hence their integral linear combinations are dense in $\R$.
Therefore there exists a closed loop $\delta$ in $M$ so that
$$0 < \, \sum_{i=1}^r \lambda_i (\delta\cdot V_i) \, < \, \lambda/3.$$
Then we may first connect the points $\Pi(w_1)$ and $\Pi(w_j)$ by an arbitrary path 
$\overline\mu_j$ in $M$ and then replace it by the path $\mu_j\, =\, k_j\delta+\overline\mu_j$
with an appropriate $k_j\in\Z$ to achieve (8.11).

The connected components $N_1, N_2, \dots, N_l$ of $N$ are in one-to-one correspondence 
with the connected components $N'_1, N'_2, \dots, N'_l$ of $N'$. Let 
$F_j: N_j' \stackrel{\subset}\longrightarrow Y\stackrel{\simeq} \longrightarrow M$
denote the composite of the inclusion and the homotopy equivalence $Y\to M$, where $j=1, \dots, l$.
We may assume that the base points $w_1, \dots, w_l\in \In(N)$ are chosen so that 
$F_j(v_j) = \Pi(w_j)$, where $v_j\in N_j'$. We may now find paths $\nu_j:[0,1]\to Y$ with
$\nu_j(0)=v_1$, $\nu_j(1)=v_j$ such that the path $F_j\circ \nu_j$ is homotopic to $\mu_j$, 
relative the endpoints.

Using the constructed 
paths $\nu_1, \nu_2, \dots, \nu_l$ in $Y$ we may define the tails of the cells of $Y$ as follows.
The tail of a cell $Q_\beta(e)\subset Y$, where $e\subset N_j'$, will start at the base point 
$v_1$, follow
the path $\nu_j$, then going along a path in $N_j'$ joining $v_j$ with an internal point of the cell
$Q_\beta(e)$, as described in section 8.8. 

It is clear that with this choice of the tails the statement of Lemma 8.9 will remain true. Indeed, in
(8.9) we will have a closed loop $\gamma_j$, having the form 
$\gamma_j = \nu_{j_1}\sigma\nu_{j_2}^{-1}$, where $\sigma$ is a path in $Y$, connecting the
base points $v_{j_1}$ and $v_{j_2}$, and such that its image in $M$ intersects only one of 
the submanifolds $V_1, \dots, V_r$, with $i=i_s$,
transversally, in the negative direction. Hence we will have $\langle\xi, [\gamma_j]\rangle <0$, 
because of (8.11).

In the next Lemma we will use the chain collapse technique (cf.  \S 7). 

\subsection{Lemma} {\it Let $\rho: \Z\pi \to \RR$ be a 
$\Sigma_\xi$-inverting ring homomorphism (cf. 1.10). Then the chain complex 
$\RR\otimes_{\Z\pi}C_\ast(\tilde Y)$ admits a sequence of $r$
collapses, so that the resulting chain complex of $\RR$-modules is free and has precisely $c_j(\omega)$ generators in each dimension $j$; these generators are in a canonical 1-1
correspondence with the index $j$ zeros of form $\omega$.
Moreover, the first $r-1$ collapses are simple.}

\begin{proof} For any proper cell $e\subset N_\alpha$ we will denote by
$i(e)\in \{1, 2, \dots, r\}$ the smalles index $i(e)=i_1$ where $\alpha =\{i_1<i_2< \dots <i_k\}$.

For any $j=1, 2, \dots, r$ we will denote by $D_j'$ the submodule of 
$\RR\otimes_{\Z\pi}C_\ast(\tilde Y)$
generated by all the basis elements $P_\beta(e)$, such that
the subset $\beta\subset \{1,2, \dots, r\}$ satisfies
$$|\beta| = r - j +1, \quad\text{and}\quad i(e)\in \beta.$$
Also, we will denote by $D_j$ the submodule of 
$\RR\otimes_{\Z\pi}C_\ast(\tilde Y)$
generated by all the basis elements $P_\beta(e)$, such that
the subset $\beta\subset \{1,2, \dots, r\}$ satisfies
$$|\beta| = r - j , \quad\text{and}\quad i(e)\notin \beta.$$
Here $e$ runs over all the cells of $\cup_{i=1}^r N_i$, cf. above. 
$D_j$ and $D_j'$ are considered as graded modules with the grading given by the dimension of the
corresponding cells, cf. (8.5).
Hence we obtain
\begin{eqnarray}
\RR\otimes_{\Z\pi}C_\ast(\tilde Y) = D_1'\oplus D_1\oplus D_2'\oplus D_2\oplus\dots
\oplus D_r'\oplus D_r\oplus F,
\end{eqnarray}
where $F$ is a graded free $\RR$-module generated by the cells of $N$ 
corresponding to the internal critical points of $g$, 
i.e. to zeros of the closed 1-form $\omega$; 
the grading of each generator of $F$ equals the index of the corresponding zero of $\omega$.
We will show that one may perform $r$ subsequent collapses in the decomposition (8.12).

First we may perform a collapse with respect to $D'_1\oplus D_1$. Indeed, from the formula
(8.10) for the boundary operator we see that only elements in $D'_1$ have their boundary
components lying in $D_1'$; in other words, we have two zeros in the first row of the matrix
(7.1). Also, the map $\gamma: D_1' \to D_1$ (in the notations of Lemma 7.1), 
i.e. the corresponding part of the differential of the chain complex $C_\ast(\tilde Y)$, is invertible.
To see this, we may label the generators $P_{\{i(e)\}\cup\beta}(e)\in D_1'$ and 
$P_{\beta}(e)\in D_1$ 
(for a given $e\subset \cup_{i=1}^r N_i$ there is at most one element of this form)
by the same letter, and then we see from (8-10) that the matrix 
representing $\gamma$ has the form $\rho(A+B)$ where $A$ is a diagonal matrix with the 
diagonal entries $\pm 1$ and $B$ is a matrix with all entries $\xi$-negative. Hence the matrix 
$\rho(A+B)$ is invertible over $\RR$ according to our assumptions about $\rho$. 

Note that, assuming that $r>1$, {\it this collapse is simple} (cf. 7.2), i.e. the map $\alpha$ 
(cf. matrix (7.1))
vanishes. Indeed, any basis element, having in the decomposition of its boundary 
a nontrivial summand in $D_1$, 
belongs to $D_1'$. Hence we will have the collapsed complex supported on the graded module
$$D_2'\oplus D_2\oplus D_3'\oplus D_3\oplus\dots
\oplus D_r'\oplus D_r\oplus F$$
and having the differential given by the same formula (8.10), in which we simply ignore the
basis elements lying in $D_1'\oplus D_1$.

Now we may perform the second collapse. The basis elements in $D_2'$ have the maximal 
length among the survived generators. Indeed, $D_2'$ contains all generators $P_\beta(e)$ with
$|\beta|=r-1$, except the ones with $\beta=\{2, 3, \dots, r\}$, which were killed during
the first collapse.
Also, the corresponding map $\gamma$ acting $D_2'\to D_2$, is again an isomorphism, because it 
can be again represented by a matrix of the form $\rho(A+B)$, where $A$ is diagonal with entries 
$\pm 1$ and where $B$ is $\xi$-negative. 

Hence, we may implement the second collapse. Let us show that, assuming 
$r>2$, the second collapse is again
simple. A basis element having in its boundary a non-trivial $D_2$-component must have the form
$P_\beta(e)$, where $|\beta|=r-1$. Such $\beta$ must contain $i(e)$ (and hence to belong to $D_2'$),
since all subsets $\beta$ with
$|\beta|=r-1$ and $i(e)\notin \beta$ were killed on the first collapse.

Suppose now that we have performed a sequence of simple collapses and have 
arrived to a chain complex
$$D_j'\oplus D_j\oplus \dots \oplus D_r'\oplus D_r\oplus F$$
with the differential given by formula (8.10), in which 
we ignore all the cells of the form $P_\beta(e)$ with
$|\beta|>r-j+1$ and also all the cells $P_\beta(e)$ with $|\beta| = r-j+1$ and $i(e)\notin \beta$, which were
killed on the previous collapses. 
Repeating the previous arguments we observe that we may perform the next collapse, killing
$D_j'\oplus D_j$. Indeed, from (8.10) we see that no element lying in 
$D_j\oplus \dots \oplus D_r'\oplus D_r\oplus F$ has a nontrivial $D_j'$-component in its boundary
(since $D_j'$ contains all survived generators with maximal $|\beta|$). In other words, we have
two zeros in the first row of matrix (7.1). The generators of $D_j'$ and $D_j$ are naturally split
into pairs (by removing index $i(e)$ from $\beta$), and if one labels the corresponding basis elements
of $D_j'$ and $D_j$ 
by the same name, we will find that the corresponding matrix $\gamma$ (cf. (7.1)) has the
form $\rho(A+B)$, where $A$ is a diagonal matrix with $\pm 1$ entries, and where $B$ is $\xi$-negative matrix. This holds since the first sum in (8.10) contains precisely one non-zero term in $D_j$.

If $j>1$ this collapse is simple, since a basis element, 
having in its boundary a non-trivial $D_j$-component must have the form
$P_\beta(e)$, where $|\beta|=r-j+1$. Such $\beta$ must contain $i(e)$, since all subsets $\beta$ with
$|\beta|=r-1$ and $i(e)\notin \beta$ were killed on the previous collapse; hence the corresponding
element $P_\beta(e)$ belongs to $D_j'$.

These arguments show that we may successfully collapse all the chain complex 
$\RR\otimes_{\Z\pi}C_\ast(\tilde Y)$ to a chain 
complex supported on graded module $F$. Note that the last collapse
$D_r'\oplus D_r\oplus F$ to $F$ {\it is not simple}, since the boundaries of the 
cells lying in the interior of $N$
(which form the basis of $F$) have components in $D_r$.

This completes the proof. \end{proof}

Lemma 8.12 clearly implies Theorem 1.11.

\bibliographystyle{amsalpha}

\enddocument